\documentclass[11pt]{article}
\usepackage[latin1]{inputenc}
\usepackage{amsmath,amsthm,amssymb}

\newtheorem{thm}{Theorem}
\numberwithin{thm}{section}

\newtheorem{prop}[thm]{Proposition}

\def\eq#1{(\ref{#1})}

\newcommand{\neweq}[1]{\begin{equation}\label{#1}}

\def\ep{\varepsilon}

\def\phi{\varphi}
\def\RR{\mathbb R}

\def\di{\displaystyle}
\def\ri{\rightarrow}
\def\intom{\int_\Omega}
\def\huo{H^1_0(\Omega )}
\def\incep{\left\{\begin{array}{cl} }
 \def\termin{\end{array}\right. }
\def\2af{2^*_\alpha}

\title{\sc Entire solutions of the nonlinear eigenvalue
 logistic problem with sign-changing potential and absorbtion}
\author{\sc Teodora-Liliana Dinu\\
\small Department of Mathematics,
``Fra\c tii Buze\c sti" College, 200585 Craiova, Romania\\ \small Email : {\tt
doruta2000@yahoo.com}}

\date{}

\begin{document}

\maketitle
\begin{abstract}
We are concerned with positive solutions decaying to zero at infinity for
 the logistic equation $-\Delta u=\lambda \left( V(x)u-f(u)\right)$
in $\RR^N$, where $V(x)$ is a variable potential that may change sign,
$\lambda$ is a real parameter, and $f$ is an
absorbtion term such that the mapping $f(t)/t$ is increasing in $(0,\infty)$.
We prove that there exists
a bifurcation non-negative number $\Lambda$ such that the above
problem has exactly one solution if $\lambda >\Lambda$,
but no such a solution exists provided $\lambda\leq\Lambda$. 
\end{abstract}

\section{Introduction and the main results}
In this paper we are concerned with the existence, uniqueness or the non-existence
of positive solutions of the
eigenvalue logistic problem with absorbtion
\neweq{log}
-\Delta u=\lambda\left(V(x)u-f(u)\right)\qquad\mbox{in}\ \RR^N,\ N\geq 3,\end{equation}
where $V$ is a smooth sign-changing potential and
 $f:[0,\infty)\ri [0,\infty)$
is a smooth function. Equations of this type arise in the study of population dynamics.
In this case,
the unknown $u$ corresponds to the density of a population,
the potential $V$ describes the birth rate
of the population, while the term $-f(u)$ in \eq{log} signifies the fact that the
population is self-limiting.
In the region where $V$ is positive (resp., negative) the population
has positive
(resp., negative) birth rate.
Since $u$ describes a population density,
we are interested in investigating only positive solutions of problem \eq{log}.

Our results are related to a certain linear eigenvalue problem. We recall in what
follows the results that we need in the sequel.
Let $\Omega$ be an arbitrary open set in $\RR^N$, $N\geq 3$.
Consider the eigenvalue problem
\neweq{eig}
 \left\{\begin{tabular}{ll}
&$-\Delta u=\lambda V(x)u$ \quad  $\mbox{in}\
\Omega\,,$\\
&$u\in\huo$. \\
\end{tabular} \right.
\end{equation}
Problems of this type have a long history. If $\Omega$ is bounded
and $V\equiv 1$, problem \eq{eig} is related to the Riesz-Fredholm theory
of self-adjoint and compact operators (see, e.g., Theorem VI.11 in \cite{b}). The case of a non-constant
potential $V$ has been first considered in the pioneering papers
of Bocher \cite{bo}, Hess and Kato \cite{hk}, Minakshisundaran and  Pleijel \cite{mp} and Pleijel \cite{p}.
For instance, Minakshisundaran and  Pleijel \cite{mp}, \cite{p} studied
the case where $\Omega$ is bounded, $V\in L^\infty (\Omega)$,
$V\geq 0$ in $\Omega$
and $V>0$ in $\Omega_0\subset\Omega$ with $|\Omega_0|>0$. An important
contribution in the study of \eq{eig}
if $\Omega$ is not necessarily bounded has been given by
Szulkin and Willem \cite{sw}
under the assumption that the sign-changing potential $V$
satisfies
$$\left\{\begin{array}{ll}
& \di V\in L^1_{\rm loc}(\Omega),\ V^+=V_1+V_2\not=0,\ V_1\in L^{N/2}(\Omega),\\
& \di \lim_{\scriptstyle x\ri y\atop\scriptstyle x\in\Omega}|x-y|^2V_2(x)=0\ \mbox{for every }y\in\overline\Omega,\
\lim_{\scriptstyle |x|\ri \infty\atop\scriptstyle
x\in\Omega}|x|^2V_2(x)=0.
\end{array}\right.\leqno (H)$$
We have denoted   $V^{+}(x)=\max\{V(x),0\}$. Obviously,
$V=V^+-V^-$, where $V^{-}(x)=\max\{-V(x),0\}$.

In order to find the principal eigenvalue of \eq{eig},
Szulkin and Willem \cite{sw} proved that
the minimization problem
$$\min\left\{\intom |\nabla u|^2dx;\ u\in \huo ,\ \intom V(x)u^2dx=1\right\}$$
has a solution $\varphi_1=\varphi_1 (\Omega)\geq 0$ which is an eigenfunction of \eq{eig}
corresponding to the eigenvalue $\lambda_1(\Omega)=\int_\Omega |\nabla\varphi_1 |^2dx$.

Throughout this paper the sign-changing potential
 $V:\RR^N\ri\RR$ is assumed to be a H\"older function that satisfies
$$
V\in L^\infty (\RR^N),\ V^+=V_1+V_2\not=0,\ V_1\in L^{N/2}(\RR^N),\ 
\lim_{|x|\ri \infty}|x|^2V_2(x)=0.
\leqno (V)$$

We suppose that the nonlinear absorbtion term $f:[0,\infty)\ri [0,\infty)$
is a $C^1$--function such that

\smallskip
\noindent $(f1)\qquad$  $f(0)=f'(0)=0$ and
$\di\liminf_{u\searrow 0}\frac{f'(u)}{u}>0$;

\smallskip
\noindent $(f2)\qquad$
the mapping $\di f(u)/u$ is increasing in $(0,+\infty)$.

\smallskip
\noindent This assumption implies
  $\lim_{u\ri +\infty}
f(u)=+\infty$. We impose that $f$ does not have a sublinear growth at infinity.
More precisely, we assume

\smallskip
\noindent $(f3)\qquad$ $\di\lim_{u\ri +\infty}
\frac{f(u)}{u}>\|V\|_{L^\infty}\,.$

\smallskip
Our framework includes the following cases: (i) $f(u)=u^2$ that corresponds to
 the Fisher equation \cite{fi}
and the  Kolmogoroff-Petrovsky-Piscounoff  equation  \cite{KPP} (see also \cite{kw}
for a comprehensive treatment of these equations); (ii) $f(u)=u^{(N+2)/(N-2)}$
(for $N\geq 6$) which is related to the conform scalar curvature equation,
cf. \cite{ln}.

For any $R>0$, denote $B_R=\{x\in\RR^N;\ |x|<R\}$ and set
\neweq{p22}\lambda_1(R)=\min\left\{\int_{B_R}|\nabla u|^2dx;\
u\in H^1_0(B_R),\ \int_{B_R}V(x)u^2dx=1\right\}\,.\end{equation}
Consequently, the mapping $R\longmapsto \lambda_1(R)$ is decreasing and so, there exists
$$\Lambda :=\lim_{R\ri\infty}\lambda_1(R)\geq 0\,.$$

We first state a sufficient condition so that $\Lambda$ is
positive. For this aim we impose the additional assumptions
\neweq{vv}
\mbox{there exists}\ \  A, \alpha >0\quad\mbox{such that}\quad
V^+(x)\leq A|x|^{-2-\alpha},\quad\mbox{ for all }\ x\in\RR^N
\end{equation}
and
\neweq{vv2003}
\lim_{x\rightarrow 0}|x|^{2(N-1)/N}V_2(x)=0.
\end{equation}

\begin{thm}\label{t1}
Assume that $V$ satisfies conditions $(V)$, \eq{vv} and \eq{vv2003}.

 Then $\Lambda>0$.
\end{thm}

Our main result asserts that $\Lambda$ plays a crucial role for the
nonlinear eigenvalue logistic problem
\neweq{p}
 \left\{\begin{tabular}{ll}
&$\di -\Delta u=\lambda \left(V(x)u-f(u)\right)$ \quad  $\mbox{in}\
\RR^N\,,$\\
&$u>0$ \quad  $\mbox{in}\ \RR^N\,,$\\
&$\di\lim_{|x|\ri\infty}u(x)=0\,.$ \\
\end{tabular} \right.
\end{equation}

The following existence and non-existence result shows that $\Lambda$
serves as a bifurcation point in our problem \eq{p}.

\begin{thm}\label{t2}
Assume that $V$ and $f$ satisfy the assumptions  $(V)$, \eq{vv},
$(f1)$, $(f2)$ and $(f3)$.

Then the following hold:

(i) problem \eq{p} has a unique solution for any $\lambda>\Lambda$;

(ii)  problem \eq{p} does not have any solution for all $\lambda\leq\Lambda$.
\end{thm}

The additional condition \eq{vv}
 implies that $V^+\in L^{N/2}(\RR^N)$, which does not follow
from the basic hypothesis $(V)$. As we shall see in the next section, this
growth assumption is essential in order to establish the existence of
positive solutions of \eq{log} {\it decaying to zero} at infinity.

In particular, Theorem \ref{t2} shows that if $V(x)<0$ for sufficiently large $|x|$
(that is, if the population has negative birth rate) then any positive
solution (that is, the population density) of \eq{log} tends to zero as $|x|\ri\infty$.

We also refer to the recent papers \cite{acm, aw, cabre, da, 
fig, ds, grossi, kyy, oss, shi, taira}
for further results related to problems of this type.

\section{Proof of Theorem \ref{t1}}
For any $R>0$, fix arbitrarily
$u\in H^1_0(B_R)$ such that
$\di\int_{ B_R}V(x)u^{2}dx=1$.
We have
$$
1=\int_{B_R}V(x)u^{2}dx
\leq\int_{B_R}V^{+}(x)u^{2}dx
=\int_{B_R}V_{1}(x)u^{2}dx+
\int_{B_R}V_{2}(x)u^{2}dx.
$$
Since $V_{1}\in L^{N/2}(\RR^N)$,
using the Cauchy-Schwarz inequality and
Sobolev embeddings we obtain
\begin{equation}\label{1.5}
\begin{array}{lll}
\di\int_{B_R}V_{1}(x)u^{2}dx\leq
\|V_{1}\|_{  L^{N/2}
(B_R)}\|u\|^{2}_{  L^{2^*}(B_R)}
\leq
C_1\| V_1\|_{  L^{{N}/{2}}(\RR^N)}
\di\int_{B_R}|\nabla u|^{2}dx,
\end{array}
\end{equation}
where $2^*=2N/(N-2)$.

Fix $\epsilon>0$.
By our assumption $(V)$, there exists positive numbers $\delta$, $R_1$ and $R$
such that
$R^{-1}<\delta <R_1<R$
 such that for all $x\in B_R$ satisfying
$|x|\geq R_1$ we have
\begin{equation}\label{1.100}
|x|^{2} V_2(x)\leq\epsilon\,.
\end{equation}
On the other hand, by $(V)$, for any $x\in B_R$ with
$|x|\leq\delta$ we have
\begin{equation}\label{1.101}
|x|^{2(N-1)/N}  V_2(x)\leq\epsilon.
\end{equation}
Define $\Omega:=\omega_1\cup\omega_2$, where
$\omega_1:=B_R\setminus\overline B_{R_1}$,
$\omega_2:=B_\delta\setminus\overline B_{1/R}$, and
$\omega:=B_{R_1}\setminus\overline B_{\delta}$.

By \eq{1.100} and Hardy's inequality we find
\begin{equation}\label{1.102}
\int_{\omega_1} V_2(x)u^2 dx\leq\epsilon\int_{\omega_1}
\frac{u^2 }{|x|^{2}}dx\leq C_2\epsilon
\int_{ B_R}|\nabla u |^2dx.
\end{equation}
Using now \eq{1.101} and H\"older's inequality we obtain
\begin{equation}\label{1.103}
\begin{array}{lll}
\di\int_{\omega_2} V_2(x)u^2 dx&\leq&\epsilon\di\int_{\omega_2}
\di\frac{u^2 }
{|x|^{2(N-1)/N}}dx\\
&\leq&\epsilon\left[\di\int_{\omega_2}
\left(\frac{1}{|x|^{2
(N-1)/N}}dx\right)^{N/2}dx\right]
^{2/N}\|u\|_{  L^{2^\star}( B_R)}^2\\
&\leq& C\epsilon\left(\di\int_{1/R}^{\delta}
\di\frac{1}{s^{N-1}}s^{N-1}\omega_N ds
\right)^{2/N}
\di\int_{ B_R}|\nabla
u |^2 dx\\
&\leq&
C_3\di\left(\delta-\frac{1}{R}\right)^{2/N}
\di\int_{ B_R}
|\nabla u |^2 dx.
\end{array}
\end{equation}
By compactness and our assumption
$(V)$, there exists a finite covering of
$\overline\omega$ by the closed balls $\overline B_{r_1}(x_1),...,
\overline B_{r_k}(x_k)$
such that, for all $1\leq j\leq k$
\begin{equation}\label{1.104}
{\rm if}\;\;|x-x_j|\leq r_j\;\; {\rm then}\;\;
|x-x_{j}|^{2(N-1)/N}  V_{2}(x)\leq\epsilon.
\end{equation}
There exists $r>0$  such that, for any $1\leq j\leq k$
$$
{\rm if}\;\; |x-x_j|\leq r\;\; {\rm then}\;\;
|x-x_j|^{2(N-1)/N} V_{2}(x)\leq\frac{\epsilon}{k}.
$$
Define $A:=\cup_{j=1}^{k}B_{r}(x_j)$. The above estimate,
 H\"older's inequality  and Sobolev embeddings yield
\begin{eqnarray*}
\int_{B_{r}(x_j)} V_2(x)u^{2} dx
&\leq&\frac{\epsilon}{k}\int_{B_{r}(x_j)}
\frac{u^{2} }{|x-x_j|^{2(N-1)/N}}dx\\
&\leq&\frac{\epsilon}{k}\left[\int_{B_{r}(x_j)}
\left(|x-x_j|^{-2(N-1)/N}
\right)^{N/2}dx
\right]^{2/
N}\|u\|^{2}_{  L^{2^{\star}}( B_R)}\\
&\leq&
C\,\frac{\epsilon}{k}\left(\int_{B_{r}}\frac{1}{|x|^{N-1}}dx\right)
^{2/N}
\int_{ B_R}|\nabla u |^2dx\\
&=&
C\,\frac{\epsilon}
{k}\left(\int_{0}^{r}\frac{1}{s^{N-1}}s^{N-1}\omega_N ds\right)
^{2/N}
\int_{ B_R}|\nabla u |^2dx\\
&=&
C' \int_{ B_R}|\nabla u |^2dx,
\end{eqnarray*}
for any $j=1,\ldots ,k$.
By addition we find
\begin{equation}\label{1.105}
\int_A  V_2(x)u^2 dx\leq
C_4\int_{ B_R}|\nabla u |^2dx.
\end{equation}
It follows from \eq{1.104} that
$ V_2\in L^{\infty}(\omega\setminus
A)$.
Actually, if $x\in\omega\setminus A$ it follows that there exists
$j\in\{1,...,k\}$ such that $r_j>|x-x_j|>r>0$.
Thus, $$ V_2(x)\leq
r^{-2(N-1)/N}\epsilon\,.$$
Hence
\begin{equation}\label{1.106}
\int_{\omega\setminus A} V_2(x)u^2 dx\leq
\epsilon
r^{-2(N-1)/N}
\int_{\omega\setminus A}u^2 dx
\leq C_5\int_{ B_R}|\nabla u |^2dx .
\end{equation}

Now from  inequalities
\eq{1.5}, \eq{1.102}, \eq{1.103}, \eq{1.105} and
\eq{1.106} we have
$$
\lambda_{1}(R)\geq
\left\{C_1\| V_1\|_{  L^{N/2}(\RR^N)}
+C_2\epsilon +C_3\left(\delta-R^{-1}\right)^{2/N}
+C_4+C_5 \right\}^{-1}
$$
and passing to the limit as $R\rightarrow\infty$ we conclude that
$$
\Lambda\geq\left(C_1\| V_1\|_{  L^{N/2}
(\RR^N)}
+C_2\epsilon +C_3\delta^{2/N}+C_4+C_5
\right)^{-1}>0.
$$
This completes the proof of Theorem \ref{t1}.
\qed

\section{An auxiliary result}
We show in this section that the logistic equation \eq{log} has entire positive solutions
if $\lambda$ is sufficiently large.
However, we are not able to establish that this solution decays to zero at infinity.
This will be proved in the next section by means of the additional assumption
\eq{vv}.
More precisely, we have

\begin{prop}\label{p1}
Assume that the functions $V$ and $f$ satisfy conditions $(V)$,
 $(f1)$,  $(f2)$ and $(f3)$. Then the problem
\neweq{finit}
 \left\{\begin{tabular}{ll}
&$\di -\Delta u=\lambda \left(V(x)u-f(u)\right)$ \quad  $\mbox{in}\
\RR^N\,,$\\
&$u>0$ \quad  $\mbox{in}\ \RR^N$\\
\end{tabular} \right.
\end{equation}
has at least one solution, for any $\lambda >\Lambda$.
\end{prop}

\begin{proof}
For any $R>0$, consider the boundary value problem
\neweq{rr}
 \left\{\begin{tabular}{ll}
&$\di -\Delta u=\lambda \left(V(x)u-f(u)\right)$ \quad  $\mbox{in}\
B_R\,,$\\
&$u>0$ \quad  $\mbox{in}\ B_R\,,$\\
&$\di u=0$ \quad  $\mbox{on}\ \partial B_R\,.$\\
\end{tabular} \right.
\end{equation}
We first prove that problem \eq{rr} has at least one solution,
for any $\lambda >\lambda_1(R)$. Indeed, the function $\overline u(x)=M$
is a supersolution of \eq{rr}, for any $M$ large enough. This follows from
$(f3)$ and the boundedness of $V$. Next, in order to find a positive subsolution, let
us consider the problem
$$\min_{u\in H^1_0(B_R)}\int_{B_R}\left(|\nabla u|^2-\lambda V(x)u^2\right)dx\,.$$
Since $\lambda >\lambda_1(R)$, it follows that the least eigenvalue $\mu_1$ is negative. Moreover, the
corresponding eigenfunction $e_1$ satisfies
\neweq{rrr}
 \left\{\begin{tabular}{ll}
&$\di -\Delta e_1-\lambda V(x)e_1=\mu_1e_1$ \quad  $\mbox{in}\
B_R\,,$\\
&$e_1>0$ \quad  $\mbox{in}\ B_R\,,$\\
&$\di e_1=0$ \quad  $\mbox{on}\ \partial B_R\,.$\\
\end{tabular} \right.
\end{equation}
Then the function $\underline u(x)=\varepsilon e_1(x)$ is a subsolution of the problem \eq{rr}. Indeed,
it is enough to check that
$$-\Delta (\varepsilon e_1)-\lambda \varepsilon Ve_1+
\lambda f(\varepsilon e_1)\leq 0\qquad
\mbox{in}\ B_R\,,$$
that is, by \eq{rrr},
\neweq{rrrr}\varepsilon\mu_1e_1+\lambda f(\varepsilon e_1)\leq 0\qquad
\mbox{in}\ B_R\,.\end{equation}
But
$$f(\varepsilon e_1)=\varepsilon f'(0)e_1+\varepsilon e_1o(1),\qquad\mbox{as}\ \varepsilon\ri 0.$$
So, since $f'(0)=0$, relation \eq{rrrr} becomes
$$\varepsilon e_1\left(\mu_1+o(1)\right)\leq 0$$
which is true, provided $\varepsilon >0$ is small enough, due to the fact that $\mu_1<0$.

Fix $\lambda >\Lambda$ and an arbitrary sequence $R_1<R_2<\ldots <R_n<\ldots$
of positive numbers such that $R_n\ri\infty$
and $\lambda_1(R_1)<\lambda$. Let $u_n$ be the
solution of \eq{rr} on $B_{R_n}$. Fix a positive number $M$ such that $f(M)/M>\|V\|_{L^\infty(\RR^N)}$.
The above arguments show that we can assume $u_{n}\leq M$ in $B_{R_n}$, for any $n\geq 1$. Since
$u_{n+1}$ is a supersolution of \eq{rr} for $R=R_n$, we can also assume that $u_n\leq u_{n+1}$ in
$B_{R_n}$. Thus the function $u(x):=\lim_{n\ri\infty} u_n(x)$ exists and is well-defined and positive in $\RR^N$.
Standard elliptic regularity arguments imply that $u$ is a solution of problem \eq{finit}.
\end{proof}

The above result shows the importance of the assumption
\eq{vv} in the statement
of Theorem \ref{t2}. Indeed, assuming that $V$ satisfies only
the hypothesis $(V)$, it is not clear whether or not the solution
constructed in the proof of Proposition \ref{p1} tends to 0 as
$|x|\ri\infty$. However, it is easy to observe that if $\lambda>\Lambda$
and $V$ satisfies \eq{vv} then problem \eq{p} has at least one solution.
Indeed, we first observe that
\neweq{subsol}\underline u(x)=\incep
&\di \ep e_1(x),\quad\mbox{if}\ x\in B_R\\
&\di 0,\quad\mbox{if}\ x\not\in B_R\termin
\end{equation}
is a subsolution of problem \eq{p}, for some fixed $R>0$, where $e_1$ satisfies
\eq{rrr}. Next, we observe that
 $\overline u(x)=n/(1+|x|^2)$ is a supersolution of \eq{p}. Indeed,
 $\overline u$ satisfies
$$-\Delta \overline u(x)=\frac{2[n(1+|x|^2)-4|x|^2]}{(1+|x|^2]^2}u(x),\qquad
x\in\RR^N.$$
It follows that $\overline u$ is a supersolution of \eq{p} provided
$$\frac{2[n(1+|x|^2)-4|x|^2]}{(1+|x|^2)^2}\geq \lambda V(x)-
\lambda f\left(\frac{n}{1+|x|^2}\right),\qquad x\in\RR^N.$$
This inequality follows from $(f3)$ and \eq{vv}, provided that
 $n$ is large enough.

\section{Proof of Theorem \ref{t2}}
We split the proof of our main result into several steps. We will assume
the conditions $(V)$, \eq{vv}, ($f1$-$f3$) are satisfied by $V$, $f$
throughout this section.

\begin{prop}\label{p2}
Let $u$ be an arbitrary solution of problem \eq{p}. Then there exists $C>0$
such that $|u(x)|\leq C |x|^{2-N}$ for all $x\in\RR^N$.
\end{prop}

\begin{proof}
Let $\omega_N$ be the surface area of the unit sphere in $\RR^N$. Consider
the function $V^+u$ as a Newtonian potential and define
$$v(x)=\frac{1}{(N-2)\omega_N}\int_{\RR^N}\frac{V^+(y)u(y)}{|x-y|^{N-2}}dy.$$
A straightforward computation shows that
\neweq{vli}-\Delta v=V^+(x)u\qquad\mbox{in}\ \, \RR^N.\end{equation}
But, by \eq{vv} and since $u$ is bounded,
$$V^+(y)u(y)\leq C |y|^{-2-\alpha},\qquad\mbox{for all}\ \, y\in\RR^N.$$
So, by Lemma~2.3 in Li and Ni \cite{ln},
$$v(x)\leq C |x|^{-\alpha},\qquad\mbox{for all }\, x\in\RR^N,$$
provided that $\alpha <N-2$.
Set $w(x)=Cv(x)-u(x)$. Hence $w(x)\ri 0$ as $|x|\ri\infty$. Let us choose $C$
sufficiently large so that $w(0)>0$. We claim that this implies
\neweq{claim} w(x)>0,\qquad\mbox{for all }\, x\in\RR^N.\end{equation}
Indeed, if not, let $x_0\in\RR^N$ be a local minimum point of $w$.
This means that $w(x_0)<0$, $\nabla w(x_0)=0$ and $\Delta w(x_0)\geq 0$. But
$$\Delta w(x_0)=-CV^+(x_0)u(x_0)+\lambda\left(V(x_0)u(x_0)-f(u(x_0))\right)< 0,$$
provided that $C>\lambda$. This contradiction implies \eq{claim}. Consequently,
$$u(x)\leq Cv(x)\leq C |x|^{-\alpha},\qquad\mbox{for any }\, x\in\RR^N.$$
So, using again \eq{vv},
$$V^+(x)u(x)\leq C |x|^{-2-2\alpha},\qquad\mbox{for all}\ \, x\in\RR^N.$$
Lemma 2.3 in \cite{ln} yields the improved estimate
$$v(x)\leq C |x|^{-2\alpha},\qquad\mbox{for all }\, x\in\RR^N,$$
provided that $2\alpha <N-2$, and so on.
Let $n_\alpha$ be the largest integer such that
$n_\alpha\alpha < N-2.$
Repeating $n_\alpha +1$ times the above argument based on Lemma~2.3 (i) and (iii)
 in \cite{ln} we obtain
$$u(x)\leq C |x|^{2-N},\qquad\mbox{for all }\, x\in\RR^N.$$
\end{proof}

\begin{prop}\label{p3}
Let $u$ be a solution of problem \eq{p}. Then $V^+u$, $V^-u$, $f(u)\in
L^1(\RR^N)$, and $u\in H^1(\RR^N)$.
\end{prop}

\begin{proof}
For any $R>0$ consider the average function
$$\overline u(R)=\frac{1}{\omega_NR^{N-1}}\int_{\partial B_R}u(x)d\sigma =
\frac{1}{\omega_N}\int_{\partial B_1}u(rx)d\sigma ,$$
where $\omega_N$ denotes the surface area of $S^{N-1}$. Then
$$\overline u'(R)=
\frac{1}{\omega_N}\int_{\partial B_1}\frac{\partial u}{\partial\nu}(rx)
d\sigma =\frac{1}{\omega_NR^{N-1}}\int_{\partial B_R}
\frac{\partial u}{\partial\nu}(x)d\sigma =\frac{1}{\omega_NR^{N-1}}\int_{B_R}
\Delta u(x)dx .$$
Hence
\begin{equation}\label{111}\begin{array}{ll}
\di\omega_NR^{N-1}\overline u'(R)=&\di -\lambda \int_{B_R}
\left( V(x)u-f(u)\right)dx=\\
&\di -\lambda \int_{B_R} V^+(x)udx+\lambda \int_{B_R}
\left( V^-(x)u+f(u)\right)dx.\end{array}
\end{equation}
 By Proposition \ref{p2},
there exists $C>0$ such that $|\overline u(r)|\leq C r^{-N+2}$, for any
$r>0$. So, by~\eq{vv},
$$\int_{1\leq |x|\leq r}V^+(x)udx\leq C A\int_{1\leq |x|\leq r}|x|^{-N-
\alpha}dx\leq C,$$
where $C$ does not depend on $r$. This implies $V^+u\in
L^1(\RR^N)$.

By contradiction, assume that $V^-u+f(u)\not\in L^1(\RR^N)$. So, by \eq{111}, $\overline u'(r)>0$
if $r$ is sufficiently large. It follows that $\overline u(r)$ does not
converge to 0 as $r\ri\infty$, which contradicts Proposition \ref{p2}. So,
$V^-u+f(u)\in L^1(\RR^N)$.
Next, in order to establish that $u\in L^2(\RR^N)$, we observe that our assumption
$(f1)$ implies the existence of some positive numbers $a$ and $\delta$ such
that $f'(t)>at$, for any $0<t<\delta$. This implies
$f(t)>at^2/2$, for any $0<t<\delta$. Since $u$ decays to 0 at infinity, it follows that the set
$\{x\in\RR^N;\ u(x)\geq\delta\}$ is compact. Hence
$$\int_{\RR^N}u^2dx=\int_{[u\geq\delta]}u^2dx+\int_{[u<\delta]}u^2dx\leq
\int_{[u\geq\delta]}u^2dx+\frac 2a\int_{[u<\delta]}f(u)dx<+\infty,$$
since $f(u)\in L^1(\RR^N)$.

It remains to prove that $\nabla u\in L^2(\RR^N)^N$.
We first observe that after multiplication by $u$ in \eq{log} and
integration we find
$$\int_{B_R} |\nabla u|^2dx-\int_{\partial B_R}u(x)
\frac{\partial u}{\partial\nu}(x)d\sigma=\lambda \int_{B_R}
\left( V(x)u-f(u)\right)dx,$$
for any $r>0$. Since $Vu-f(u)\in L^1(\RR^N)$, it follows that the left
hand-side has a finite limit as $r\ri\infty$.
Arguing by contradiction and assuming that $\nabla u\not\in L^2(\RR^N)^N$, it
follows that there exists $R_0>0$ such that
\neweq{112}\int_{\partial B_R}u(x)
\frac{\partial u}{\partial\nu}(x)d\sigma\geq \frac 12\int_{B_R} |\nabla u|^2dx,
\qquad\mbox{for any $R\geq R_0$}.
\end{equation}

Define the functions
$$A(R)=\int_{\partial B_R}u(x)
\frac{\partial u}{\partial\nu}(x)d\sigma,\qquad B(R)=
\int_{\partial B_R} u^2(x)d\sigma,\qquad
C(R)=\int_{B_R} |\nabla u(x)|^2dx.$$
Relation \eq{112} can be rewritten as
\neweq{11112}
A(R)\geq \frac 12 \, C(R),\qquad\mbox{for any $R\geq R_0$}.
\end{equation}
On the other hand, by the Cauchy-Schwarz inequality,
$$A^2(R)\leq \left(\int_{\partial B_R} u^2d\sigma\right)
\left(\int_{\partial B_R} \left|\frac{\partial u}{\partial\nu}\right|^2
d\sigma\right)\leq B(R)C'(R).$$
Using now \eq{11112} we obtain
$$C'(R)\geq\frac{C^2(R)}{4B(R)},\qquad\mbox{for any }\, R\geq R_0.$$
Hence
\neweq{113}\frac{{\rm d}}{{\rm d}r}\left[\frac{4}{C(r)}+\int_0^r\frac{dt}{B(t)}
\right]_{r=R}\leq 0,\qquad\mbox{for any }\, R\geq R_0.\end{equation}
But, since $u\in L^2(\RR^N)$, it follows that $\int_0^\infty B(t)dt$ converges, so
\neweq{114}\lim_{R\ri\infty}\int_0^R\frac{dt}{B(t)}=+\infty.
\end{equation}
On the other hand, our assumption $|\nabla u|\not\in L^2(\RR^N)$ implies
\neweq{115}\lim_{R\ri\infty}\frac{1}{C(R)}=0.
\end{equation}
Relations \eq{113}, \eq{114} and \eq{115} yield a contradiction, so our proof
is complete.
\end{proof}

\begin{prop}\label{p4}
Let $u$ and $v$ be two distinct solutions of problem \eq{p}. Then
$$\lim_{R\ri\infty}\int_{\partial B_R}u(x)\frac{\partial v}{\partial\nu}(x)
d\sigma =0.$$
\end{prop}

\begin{proof}
By multiplication with $v$ in \eq{p} and integration on $B_R$ we find
$$\int_{B_R}\nabla u\cdot\nabla vdx-\int_{\partial B_R}u
\frac{\partial v}{\partial\nu}
d\sigma= \lambda \int_{B_R}
\left( V(x)uv-f(u)v\right)dx.$$
So, by Proposition \ref{p3}, there exists and is finite $\lim_{R\ri\infty}\int_{\partial B_R}
u\frac{\partial v}{\partial\nu}d\sigma$. But, by the Cauchy-Schwarz inequality,
\neweq{cs}\left|\int_{\partial B_R}
u\frac{\partial v}{\partial\nu}d\sigma\right|\leq\left(\int_{\partial B_R}
u^2d\sigma\right)^{1/2}\left(\int_{\partial B_R}
|\nabla v|^2d\sigma\right)^{1/2}.
\end{equation}
Since $u$, $|\nabla v|\in L^2(\RR^N)$, it follows that
$\int_0^\infty\left(\int_{\partial B_R}(u^2+|\nabla v|^2)d\sigma\right)dx$ is
convergent. Hence
\neweq{cs2}\lim_{R\ri\infty}\int_{\partial B_{R}}(u^2+|\nabla v|^2)d\sigma
=0.\end{equation}
Our conclusion now follows by \eq{cs} and \eq{cs2}.
\end{proof}

\medskip
{\sc Proof of Theorem \ref{t2}}. (i) The existence of a solution follows with the arguments given in the preceding section.
In order to establish the uniqueness, let $u$ and $v$ be two solutions of \eq{p}.
We can assume without loss of generality that $u\leq v$. This follows from the
fact that $\overline u=\min\{u,v\}$ is a supersolution of \eq{p} and $\underline u$
defined in \eq{subsol} is an arbitrary small subsolution. So, it sufficient
to consider the ordered pair consisting of the corresponding solution and $v$.

Since $u$ and $v$ are solutions we have, by Green's formula,
$$\int_{\partial B_R}\left( u\frac{\partial v}{\partial\nu}
 -v\frac{\partial u}{\partial\nu}\right)d\sigma=\lambda
 \int_{B_R}uv\left(\frac{f(v)}{v}-\frac{f(u)}{u}\right)dx.$$
 By Proposition \ref{p4}, the left hand-side converges to 0 as $R\ri\infty$.
 So,  $(f1)$ and our assumption $u\leq v$ force $u=v$ in $\RR^N$.

(ii) By contradiction, let $\lambda\leq\Lambda$ be such that problem \eq{p} has a solution
for this $\lambda$. So
$$\int_{B_R}|\nabla u|^2dx-\int_{\partial B_R}u\frac{\partial u}{\partial\nu}d\sigma
=\lambda \int_{B_R}
\left( V(x)u^2-f(u)u\right)dx.$$
By Propositions \ref{p3} and \ref{p4} and letting $R\ri\infty$ we find
\neweq{vic}\int_{\RR^N}|\nabla u|^2dx<
\lambda \int_{\RR^N} V(x)u^2dx.\end{equation}
On the other hand, using the definition of $\Lambda$ and \eq{p22} we obtain
\neweq{vicv}\Lambda\int_{\RR^N} V\zeta^2dx\leq \int_{\RR^N} |\nabla \zeta|^2dx,
\end{equation}
for any $\zeta\in C^2_0(\RR^N)$ such that $\int_{\RR^N} V\zeta^2dx>0$.

Fix $\zeta\in C^2_0(\RR^N)$ such that $0\leq\zeta\leq 1$, $\zeta (x)=1$
if $|x|\leq 1$, and $\zeta (x)=0$
if $|x|\geq 2$. For any $n\geq 1$ define $\Psi_n(x)=\zeta_n(x)u(x)$,
where $\zeta_n(x)=\zeta (|x|/n)$. Thus
$\Psi_n(x)\ri u(x)$ as $n\ri\infty$, for any $x\in\RR^N$. Since $u\in H^1(\RR^N)$,
it follows by Corolarry IX.13 in \cite{b}
that $u\in L^{2N/(N-2)}(\RR^N)$. So, the Lebesgue dominated convergence theorem
yields
$$\Psi_n\ri u\qquad\mbox{in }\, L^{2N/(N-2)}(\RR^N).$$
We claim that
\neweq{cc}\nabla\Psi_n\ri\nabla u\qquad\mbox{in }\, L^{2}(\RR^N)^N.
\end{equation}
Indeed, let $\Omega_n:=\{x\in\RR^N;\ n<|x|<2n\}$. Applying H\"older's
inequality we find
\neweq{co}\begin{array}{ll}\di
\|\nabla \Psi_n-\nabla u\|_{L^2(\RR^N)}\leq&\di
\|(\zeta_n-1)\nabla u\|_{L^2(\RR^N)}+\|u\nabla\zeta_n\|_{L^2(\Omega_n)}\leq\\
&\di\|(\zeta_n-1)\nabla u\|_{L^2(\RR^N)}+\|u\|_{L^{2N/(N-2)}(\Omega_n)}\cdot
\|\nabla\zeta_n\|_{L^{N}(\RR^N)}.\end{array}
\end{equation}
But, since $|\nabla u|\in L^2(\RR^N)$, it follows by Lebesgue's dominated convergence
theorem that
\neweq{co11}\lim_{n\ri\infty}\|(\zeta_n-1)\nabla u\|_{L^2(\RR^N)}=0.
\end{equation}
Next, we observe that
\neweq{co12}\|\nabla\zeta_n\|_{L^{N}(\RR^N)}=
\|\nabla\zeta\|_{L^{N}(\RR^N)}.
\end{equation}
Since $u\in L^{2N/(N-2)}(\RR^N)$ then
\neweq{co13}\lim_{n\ri\infty}\|u\|_{L^{2N/(N-2)}(\Omega_n)}=0.\end{equation}
Relations \eq{co}--\eq{co13} imply our claim \eq{cc}.

Since $V^{\pm}u^2\in L^1(\RR^N)$ and $V^{\pm}\Psi_n^2\leq V^{\pm}u^2$, it
follows by Lebesgue's dominated convergence theorem that
$$\lim_{n\ri\infty}\int_{\RR^N}V^{\pm}\Psi_n^2dx=\int_{\RR^N}V^{\pm}u^2dx.$$
Consequently
\neweq{dor}\lim_{n\ri\infty}\int_{\RR^N}V\Psi_n^2dx=\int_{\RR^N}Vu^2dx.
\end{equation}
So, by \eq{vic} and \eq{dor}, it follows that there exists $n_0\geq 1$ such that
$$\int_{\RR^N}V\Psi_n^2dx>0,\qquad\mbox{for any }\, n\geq n_0.$$
This means that we can write \eq{vicv} for $\zeta$ replaced by
$\Psi_n\in C^2_0(\RR^N)$. Using then \eq{cc} and \eq{dor}
 we find
\neweq{final}\int_{\RR^N}|\nabla u|^2dx\geq\Lambda
\int_{\RR^N}Vu^2dx.
\end{equation}
Relations \eq{vic} and \eq{final} yield a contradiction, so problem \eq{p} has no
solution if $\lambda\leq\Lambda$. \qed

\medskip
{\bf Acknowledgments}. The author is grateful to Professor
Congming Li for his interest in this work and for numerous
comments and suggestions on a first version of this paper.

\end{document}